\documentclass[12pt]{article}

\usepackage{amssymb,amsmath}
\usepackage{graphicx}
\usepackage{tikz-cd}

\usepackage{pstricks}

\usepackage[utf8]{inputenc}

\title{\bf\Large An Embedding of Schwartz Distributions in the Algebra of Asymptotic Functions}
\author{Michael Oberguggenberger\\[5pt]
Institut f\"{u}r Mathematik und Geometrie\\
Universitat Innsbruck\\
A-6020 Inssbruck, Austria\\
michael@mat1.uibk.ac.at\\
and\\[5pt]
Todor Todorov\\[5pt]
Mathematics Department\\
California Polytechnic State University\\
San Luis Obispo, California 93407, USA\\
ttodorov@calpoly.edu\\[5pt]
(Received October 28, 1996 and in revised form June 21, 1997)
}

\date{}
\begin{document}

\maketitle

\begin{abstract} We present a solution of the problem of multiplication of Schwartz distributions by
embedding the space of distributions into a differential algebra of generalized functions, called in the
paper ``asymptotic function'', similar to but different from J. F. Colombeau's algebras of new generalized
functions.
\end{abstract}


\medskip


\section{Introduction}

The main purpose of this paper is to prove the existence of an embedding $\Sigma_{{D},\Omega}$ of the space of
Schwartz distributions ${\cal D}^\prime(\Omega)$ into the algebra of asymptotic functions $^\rho \mathcal E(\Omega)$ which preserves all linear
operations in ${\cal D}^\prime(\Omega)$. Thus, we offer a solution of the problem of multiplication of Schwartz distributions
since the multiplication within ${\cal D}^\prime(\Omega)$ is impossible (L. Schwartz [1]).

	The algebra $^\rho \mathcal E(\Omega)$ is defined in the paper as a factor space of nonstandard smooth functions. The field of the scalars $^\rho\mathbb C$ of the algebra $^\rho\mathcal E(\Omega)$ coincides with the complex counterpart of the Robinson asymptotic numbers -known also as {\it Robinson's field with valuation} (Robinson's
(2]) and Lightstone and Robinson [3]). The embedding $\Sigma_{D,\Omega}$ is constructed in the form
$\Sigma_{D,\Omega}=Q_\Omega \circ D\star\Pi\bullet\, ^*$,  where (in backward order): $^*$ is the extension mapping (in the sense of nonstandard analysis), $\bullet$ is the Schwartz multiplication in ${\cal D}^\prime(\Omega)$ (more precisely, its nonstandard
extension), $\star$ is the convolution operator (more precisely, its nonstandard extension), $\circ$
denotes ``composition'', $Q_\Omega$ is the quotient mapping (in the definition of the algebra of asymptotic
functions) and $D$ and $\Pi_\Omega$ are fixed nonstandard internal functions with special properties whose existence
is proved in this paper.

Our interest in the algebra $^\rho \mathcal E(\Omega)$ and the embedding ${\cal D}^\prime(\Omega)\subset\, ^\rho \mathcal E(\Omega)$, is due to their role in the problem of multiplication of Schwartz distributions, the nonlinear theory of generalized functions and its
applications to partial differential equations (Oberguggenberger [4]), (Todorov [5] and [6]). In
particular, there is a strong similarity between the algebra of asymptotic functions $^\rho \mathcal E(\Omega)$ and its
generalized scalars $^\rho\mathbb{C}$, discussed in this paper, and the algebra of generalized functions ${\cal G}(\Omega)$ and their
generalized scalars $\overline{\mathbb{C}}$, introduced by J. F. Colombeau in the framework of standard analysis
(Colombeau [7], pp. 63, 138 and Colombeau [8], \S 8.3, pp.161-166). We should mention that
the involvement of nonstandard analysis has resulted in some improvements of the corresponding
standard counterparts; one of them is that $^\rho\mathbb{C}$ is an algebraically closed field while its standard
counterpart $\overline{\mathbb{C}}$ in J. F. Colombeau's theory is a ring with zero divisors.

This paper is a generalization of some results in [9] and [10] (by the authors of this paper,
respectively) where only the embedding of the tempered distributions $\mathcal S^\prime(\mathbb{R}^d)$ in $^\rho \mathcal E(\mathbb{R}^d)$ has been
established. The embedding of all distributions ${\cal D}^\prime(\Omega)$, discussed in this paper, presents an essentially
different situation. We should mention that the algebra $^\rho \mathcal E(\mathbb{R}^d)$ was recently studied in (Hoskins and Sousa Pinto [11]).

Here $\Omega$ denotes an open set of $\mathbb{R}^d$ ($d$ is a natural number), $\mathcal E(\Omega) =C^\infty(\Omega)$
and ${\cal D}(\Omega)=C^\infty_0(\Omega)$ denote the usual classes of $C^\infty$-functions on $\Omega$ and
$C^\infty$-functions with compact support in $\Omega$ and ${\cal D}^\prime(\Omega)$, and $\mathcal E^{\prime}(\Omega)$ denote the
classes of Schwartz distributions on $\Omega$ and Schwartz distributions with compact
support in $\Omega$, respectively. As usual, $\mathbb{N}$, $\mathbb{R}$, $\mathbb{R}_+$, and $\mathbb{C}$ will be the systems of the natural, real, positive real and complex numbers, respectively, and we use also the notation $\mathbb N_0 = \{0\}\cup \mathbb{N}$. For the partial
derivatives we write $\partial^\alpha$, $\alpha\in \mathbb{N}_0^d$. If $\alpha = (\alpha_1,\dots , \alpha_d)$ for some $\alpha\in \mathbb{N}_0^d$, we write $|\alpha|=\alpha_1+\cdots +\alpha_d$
and if $x= (x_1,\dots x_d)$ for some $x\in\mathbb{R}^d$, we write $x^\alpha=x_1^{\alpha_1}x_2^{\alpha_2}\dots x_d^{\alpha_2}\dots x_d^{\alpha_d}$. For a general reference to distribution theory we refer to Bremermann [12] and Vladimirov [13].

Our framework is a nonstandard model of the complex numbers $\mathbb{C}$, with degree of saturation larger
than card$(\mathbb{N})$. We denote by $^*\mathbb{R}$, $^*\mathbb{R}_+$, $^*\mathbb{C}$,
$^* \mathcal E(\Omega)$ and $^*{\cal D}(\Omega)$ the nonstandard extensions of $\mathbb{R}$, $\mathbb{R}_+$, $\mathbb{C}$,
$\mathcal E(\Omega)$ and ${\cal D}(\Omega)$, respectively. If $X$ is a set of complex numbers or a set of (standard) functions, then
$^*X$ will be its nonstandard extension and if $f: X\to Y$ is a (standard) mapping, then ${^*f}\colon {^*X}\to
\, ^\star Y$ will be its nonstandard extension. For integration in $^*\mathbb{R}^d$ we use the $*$-Lebesgue integral.
We shall often use the same notation, $\| x\|$, for the Euclidean norm in $\mathbb{R}^d$ and its nonstandard extension in
$^\star\mathbb{R}^d$. For a short introduction to nonstandard analysis we refer to the Appendix in Todorov [6]. For a more detailed
exposition we recommend Lindstrom [14], where the reader will find many references to the subject.

\section{Test Functions and Their Moments}

	\quad\; In this section we study some properties of the test functions in ${\cal D}(\mathbb{R}^d)$ (in a standard setting) which
we shall use subsequently.

Following (Colombeau [7], p. 55), for any $k\in\mathbb{N}$ we define the set of test functions:
$$
A_k=\{\varphi\in {\cal D}(\mathbb{R}^d)\colon \varphi \mbox{ is real-valued, } \varphi(x) = 0 \mbox{ for } \| x \| > 1;
$$
$$
\int_{\mathbb{R}^d}\varphi(x)dx=1 \mbox{ and } \int_{\mathbb{R}^d}x^\alpha\varphi(x)dx=0 \mbox{ for }\alpha\in \mathbb{N}_0^d, 1\leq |\alpha |\leq k\}.\leqno{(2.1)}
$$
Obviously, $A_1\supset A_2\supset A_3\supset\dots$. Also, we have $A_k\neq\emptyset$ for all $k\in\mathbb{N}$ (Colombeau [7], Lemma
(3.3.1), p. 55).

In addition to the above we have the following result:

\medskip

{\bf Lemma 2.2.} For any $k\in\mathbb{N}$
$$
\inf_{\varphi\in A_k}\left(\int_{\mathbb{R}^d}|\varphi(x)|dx\right)=1. \leqno{(2.2)}
$$
More precisely, for any positive real $\delta$ there exists $\varphi$ in $A_k$ such that
$$
1\leq \int_{\mathbb{R}^d}|\varphi(x)|dx\leq 1+\delta.
$$
In addition, $\varphi$ can be chosen symmetric.

\medskip

{\bf {\bf Proof.}} We consider the one dimensional case $d = 1$ first. Start with some fixed positive (real
valued) $\psi$ in ${\cal D}(\mathbb{R})$ such that $\psi(x) = 0$ for $|x|\geq 1$ and $\int_{\mathbb{R}}\psi(x)dx = 1$ ($\psi$ can be also chosen symmetric if needed). We shall look for $\varphi$ in the form:
$$
\varphi(x)=\sum_{j=0}^{k}c_j\psi\left(\frac{x}{\varepsilon^j}\right),
$$
$x\in \mathbb{R}$, $\varepsilon\in\mathbb{R}_+$. We have to find $c_j$ for which $\varphi\in A_k$. Observing that
$$
\int_{\mathbb{R}}x^j\psi\left(\frac{x}{\varepsilon^j}\right)dx=\varepsilon^{(i+1)j}\int_{\mathbb{R}}y^j\psi(y)dy,
$$
for $i = 0,1,\dots ,k$, we derive the system for linear equations for $c_j$:
$$
\left\{\begin{array}{ll}
\displaystyle \sum_{j=0}^{k}c_j\,\varepsilon^j=1,\\[8pt]
\displaystyle \left(\int_{\mathbb{R}}y^i\psi(y)dy\right)\sum_{j=0}^{k}c_j\,\varepsilon^{(i+1)j}=0,\; i=1,\dots ,k.
\end{array}\right.
$$
The system is certainly satisfied, if
$$
\left\{\begin{array}{ll}
\displaystyle\sum_{j=0}^{k}c_j\,\varepsilon^j=1,\\[8pt]
\displaystyle\sum_{j=0}^{k}c_j\,\varepsilon^{(i+1)j}=0, i=1,\dots ,k,
\end{array}\right.
$$
which can be written in the matrix form $V_{k+1}(\varepsilon)\,C=B$, where $V_{k+1}(\varepsilon)$ is Vandermonde
$(k+1)\times (k+1)$ matrix, $C$ is the column of the unknowns $c_j$, and $B$ is the column whose top entry is 1 and all others are 0.
For the determinant we have det $V_{k+1}(\,\varepsilon)\neq 0$ for $\,\varepsilon\neq 1$, therefore, the system has a unique solution
$(c_0,c_1,c_2,\dots ,c_k)$. Our next goal is to show that this solution is of the form:
$$
c_j=\pm\frac{\,\varepsilon^{\alpha_j}(1+\,\varepsilon P_j(\varepsilon))}{\,\varepsilon^{\beta}(1+\,\varepsilon P(\varepsilon))},\leqno{(2.3)}
$$
where $P_j$ and $P$ are polynomials and
$$
\alpha_j=\sum_{q=1}^{k-1}q(k+1-q)+\sum_{m=j}^{k-1}(k+1-m), \leqno{(2.4)}
$$
for $0\leq j\leq k$, and
$$
\beta=k+\sum_{q=1}^{k-1}q(k+1-q).
$$
The coefficients $c_0,c_1,c_2,\dots ,c_k$ will be found by Cramer's rule. The formula for Vandermonde determinants
gives
\begin{eqnarray*}
\displaystyle\Pi_{m=1}^{k}\Pi_{q=m+1}^{k+1}(\varepsilon^q-\varepsilon^m)&=&\Pi_{m=1}^{k}\left(\Pi_{q=m+1}^{k+1}\varepsilon^m(\varepsilon^{q-m}-1)\right)\\[8pt]
&=&\varepsilon^\beta\Pi_{m=1}^{k}\Pi_{q=m+1}^{k+1}(\varepsilon^{q-m}-1)=\pm\varepsilon^\beta(1+\varepsilon P(\varepsilon)),
\end{eqnarray*}
for some polynomial $P$, where
$$
\beta=\sum_{m=1}^{k}m(k+1-m)=k+\sum_{m=1}^{k-1}m(k+1-m).
$$
To calculate the numerator in (2.3), we have to replace the $j$th column of the matrix by the column $B$
(whose top entry is 1 and all others are 0) and calculate the resulting determinant $D_j$. Consider first the
case $1\leq j\leq k-1$. By developing with respect to the $j$th column, we get
$$
D_j=\pm\det\left(
\begin{array}{ccccccc}
1,&\varepsilon^2,&\dots & \varepsilon^{2(j-1)},&\varepsilon^{2(j+1)},&\dots &\varepsilon^{2k}\\
1,&\varepsilon^3,&\dots & \varepsilon^{3(j-1)},&\varepsilon^{3(j+1)},&\dots &\varepsilon^{3k}\\
\dots &\dots &\dots & \dots &\dots &\dots &\dots \\
1,&\varepsilon^{k+1},&\dots & \varepsilon^{(k+1)(j-1)},&\varepsilon^{(k+1)(j+1)},&\dots &\varepsilon^{(k+1)k}
\end{array}
\right).
$$
We factor out $\varepsilon^2,\varepsilon^4,\dots, \varepsilon^{2(j-1)},\varepsilon^{2(j+1)},\dots,\varepsilon^{2k}$ and obtain:
\begin{eqnarray*}
D_j&=&\pm \varepsilon^{1\cdot 1}\varepsilon^{2\cdot 2}\dots \varepsilon^{(j-1)\cdot 2}\varepsilon^{(j+1)\cdot 2}\dots \varepsilon^{2k}\\
&\times\det&\left(
\begin{array}{cccccccc}
1,&1,&1,& \dots &1,&1, &\dots &1\\
1,&\varepsilon,& \varepsilon^2&\dots & \varepsilon^{(j-1)},&\varepsilon^{(j+1)},&\dots &\varepsilon^{k}\\
\dots &\dots &\dots & \dots &\dots &\dots &\dots \\
1,&\varepsilon^{k-1},&\varepsilon^{2(k-1)},  &\dots & \varepsilon^{(j-1)(k-1)},&\,\varepsilon^{(j+1)(k-1)},&\dots &\varepsilon^{k(k-1)}
\end{array}
\right).
\end{eqnarray*}
The latter is a Vandermonde determinant again, and we have
\begin{eqnarray*}
D_j=&\pm&\varepsilon^{1\cdot 1+2\cdot 2+\cdots+ (j-1)2+(j+1)2+\cdots +k\cdot 2}\\
&\times&
(\varepsilon-1)(\varepsilon^2-1)\dots (\varepsilon^{j-1}-1)(\varepsilon^{j+1}-1)\dots (\varepsilon^k-1)\\
&\times& (\varepsilon^2-\varepsilon)(\varepsilon^3-\varepsilon)\dots (\varepsilon^{j-1}-\varepsilon)(\varepsilon^{j+1}-\varepsilon)\dots (\varepsilon^k-\varepsilon)\\
&\times& \dots\dots\dots\dots\dots\dots\dots\dots \\
&\times& (\varepsilon^{j-1}-\varepsilon^{j-2})(\varepsilon^{j+1}-\varepsilon^{j-2})\dots (\varepsilon^{k}-\varepsilon^{j-2})\\
&\times& (\varepsilon^{j+1}-\varepsilon^{j-1})\dots (\varepsilon^{k}-\varepsilon^{j-1})\\
&\times& \dots\dots\dots\dots \\
&\times& (\varepsilon^{k}-\varepsilon^{k-1}).
\end{eqnarray*}
Hence, factoring out $\varepsilon^{(i-1)(k-i)}$ in the $i$th row above, we get $D_j=\pm\varepsilon^{\alpha_j}(1+\varepsilon P_j(\varepsilon))$ for some
polynomials $P_j(\varepsilon)$ and
\begin{eqnarray*}
\alpha_j&=& 1\cdot 1+2\cdot 2+\cdots+ (j-1)\cdot 2+(j+1)\cdot 2+\cdots +k\cdot 2\\
&+&1\cdot(k-2)+2\cdot(k-3)+\cdots +(j-1)\cdot(k-j)\\
&+&(j+1)\cdot(k-j-1)+\cdots + (k-1)\cdot 1\\
&=&1\cdot k+2\cdot (k-1)+\cdots +(j-1)\cdot (k-j+2)+(j+1)\cdot (k-j+1)+\cdots +(k-1)\cdot 3+k\cdot 2\\
&=&\sum_{q=1}^{j-1}q(k+1-q)+\sum_{m=j}^{k-1}(m+1)(k+1-m)=\sum_{q=1}^{k-1}q(k+1-q)+\sum_{m=j}^{k-1}(k+1-m),
\end{eqnarray*}
which coincides with the desired result (2.4) for $\alpha_j$, in the case $1\leq j\leq k-1$. For the extreme cases
$j =0$ and $j = k$, we obtain
$$
\alpha_0=\sum_{m=0}^{k-1}(m+1)(k+1-m)=\sum_{q=1}^{k-1}q(k+1-q)+\sum_{m=0}^{k-1}(k+1-m),
$$
$$
\alpha_k=\sum_{q=1}^{k-1}q(k+1-q),
$$
which both can be incorporated in the formula (2.4) for $\alpha_j$. Finally, Cramer's rule gives the expression
(2.3) for $c_j$.

Now, taking into account that $\psi\geq 0$, by assumption, and the fact that $|1+\varepsilon P(\varepsilon)| > |1-|\varepsilon P(\varepsilon)|| =
1-\varepsilon |P(\varepsilon)| > 0$ for all sufficiently small epsilon, we obtain
$$
\int_{\mathbb{R}}|\varphi(x)|dx\leq\sum_{j=0}^{k}|c_j|\,\varepsilon^j\leq\sum_{j=0}^{k}
\frac{\varepsilon^{j+\alpha_j}(1+\varepsilon|P_j(\varepsilon)|)}{\varepsilon^\beta(1-\varepsilon |P(\varepsilon)|)},
$$
and this latter expression can be made smaller than $1+\delta$ for sufficiently small $\varepsilon$ if a) $j +\alpha_j-\beta>0$ for
$0\leq j\leq k-1$, and b) $k+\alpha_k-\beta=0$. Now, b) is obvious, as for a), we have:
$$
j+\alpha_j-\beta=j+\sum_{m=j}^{k-1}(k+1-m)-k=\frac{1}{2}(k-j)(k-j+1)>0,
$$
for $0\leq j\leq k-1$. To generalize the result for arbitrary dimension $d$, it suffices to consider a product of
functions of one real variable. The proof is complete. $\square$

\section{Nonstandard Delta Functions}

\quad\; We prove the existence of a nonstandard function $D$ in $^*{\cal D}(\mathbb{R}^d)$ with special properties. The proof is
based on the result of Lemma 2.2 and the Saturation Principle (Todorov [6], p. 687). We also
consider a type of nonstandard cut-off-functions which have close counterparts in standard analysis. The
applications of these functions are left for the next sections.

\medskip

{\bf Lemma 3.1} (Nonstandard Mollifiers). For any positive infinitesimal $\rho$ in $^*\mathbb{R}$ there exists a
nonstandard function $\theta$ in $^\star{\cal D}(\mathbb{R}^d)$ with values in $^*\mathbb{R}$, which is symmetric and which satisfies the following properties:

\smallskip

(i) $\theta(x)=0$ for $x\in\, ^*\mathbb{R}^d, \|x\|>1$,

\smallskip

(ii) $\int_{^\star\mathbb{R}^d}\theta(x)dx=1$,

\smallskip

(iii) $\int_{^\star\mathbb{R}^d}\theta(x)x^\alpha dx=0$, for all $\alpha\in \mathbb{N}_0^d$, $\alpha\neq 0$,

\smallskip

(iv) $\int_{^\star\mathbb{R}^d}|\theta(x)|dx\approx 1$,

\smallskip

(v) $\displaystyle|\ln \rho|^{-1}\left(\sup_{x\in\, ^\star\mathbb{R}^d}|\partial^\alpha \theta(x)|\right)\approx 0$ for all a $\alpha\in \mathbb{N}_0^d$,

\smallskip
\noindent where $\approx$ is the infinitesimal relation in $^*\mathbb{C}$. We shall call this type of function {\it nonstandard $\rho$-mollifiers}

\medskip

{\bf Proof.} For any $k \in\mathbb{N}$, we define the set of test functions:
$$
\overline{A}_k = \{\varphi\in{\cal D}(\mathbb{R}^d)\colon \varphi\mbox{ is real-valued and symmetric},
$$
$$
\varphi(x)=0 \mbox{ for } \|x\|\geq 1, \int_{\mathbb{R}^d}\varphi(x)dx=1,
$$
$$
\int_{\mathbb{R}^d}x^\alpha\varphi(x)dx=0 \mbox{ for } 1\leq |\alpha|\leq k, \int_{\mathbb{R}^d}|\varphi(x)|dx<1+\frac{1}{k}\}
$$
and the internal subsets of $^*{\cal D}(\mathbb{R}^d)$:
$$
{\cal A}_k=\{\varphi\in\, ^*(\overline{A}_k):|\ln\, \rho|^{-1}\left( \sup_{x\in\, ^*\mathbb{R}^d} |\partial^\alpha(^*\!\varphi(x))|\right)                    <\frac{1}{k} \mbox{ for } |\alpha|\leq k\}.
$$
Obviously, we have $\overline{A}_1\supset \overline{A}_2\supset \overline{A}_3\supset\dots $ and ${\cal A}_1\supset {\cal A}_2\supset {\cal A}_3\supset\dots $.
Also we have $\overline{A}_k\neq 0$ for all $k\in\mathbb{N}$ by Lemma 2.2.
On the other hand, we have $\overline{A}_k\subset {\cal A}_k$ in the sense that $\varphi \in \overline{A}_k$ implies
$^*\!\varphi \in {\cal A}_k$, since
$$
\sup_{x\in\, ^\star\mathbb{R}^d} |\partial^\alpha(^*\!\varphi(x))| = \sup_{x\in\, \mathbb{R}^d}|\partial^\alpha\varphi(x)|=
\sup_{x\leq 1} |\partial^\alpha(^*\!\varphi(x))|,
$$
is a real (standard) number and, hence, $\displaystyle|\ln\, \rho|^{-1}\left( \sup_{x\in\,^\star\mathbb{R}^d} |\partial^\alpha(^*\!\varphi(x))|\right)$ is infinitesimal. Thus, we have ${\cal A}_k\neq\emptyset$
for all $k$ in $\mathbb{N}$. By the Saturation Principle (Todorov [6], p. 687), the intersection
$\displaystyle{\cal A}=\bigcap_{k\in \mathbb{N}}{\cal A}_k$ is non-empty and thus, any $\theta$ in ${\cal A}$ has the desired properties.\ \ $\square$

\medskip

{\bf Definition 3.2} ($\rho$-Delta Function). Let $\rho$ be a positive infinitesimal. A nonstandard function $D$
in $^* {\cal D}(\mathbb{R}^d)$ is called a $\rho$-{\it delta function} if it takes values in $^* \mathbb{R}$, it is symmetric and it satisfies the following properties:

\smallskip

(i) $D(x)=0$ for $x\in\, ^*\mathbb{R}^d, \|x\|\geq \rho$;

\smallskip

(ii) $\int_{^\star\mathbb{R}^d}D(x)\,dx=1$,

\smallskip

(iii) $\int_{^\star\mathbb{R}^d}D(x)x^\alpha \,dx=0$, for all $\alpha\in \mathbb{N}_0^d$, $\alpha\neq 0$,

\smallskip

(iv) $\int_{^\star\mathbb{R}^d}|D(x)|\,dx\approx 1$,

\smallskip

(v) $\displaystyle|\ln \rho|^{-1}\left(\rho^{d+|\alpha|}\sup_{x\in{^*\mathbb{R}^d}}|\partial^\alpha D(x)|\right)\approx 0$ for all $\alpha\in \mathbb{N}_0^d$;

\medskip

{\bf Theorem 3.3} (Existence). For any positive infinitesimal $\rho$ in $^*\mathbb{R}$ there exists a $\rho$-delta function.

\medskip

{\bf Proof.} Let $\theta$ be a nonstandard $\rho$-mollifier of the type described in Lemma 3.1. Then the
nonstandard function $D$ in $^*{\cal D}(\mathbb{R}^d)$, defined by
$$
D(x) = \rho^{-d}\,\theta(x/\rho),\ \ \  x\in\, ^*\mathbb{R}^d \leqno{(3.1)}
$$
satisfies (i)-(v). $\square$

\medskip

{\bf Remark.} The existence of nonstandard functions $D$ in $^*{\cal D}(\mathbb{R}^d)$ with the above properties
is in sharp contrast with the situation in standard analysis where there is no $D$ in ${\cal D}(\mathbb{R}^d)$
which satisfies both (ii) and (iii). Indeed, if we assume that $D$ in ${\cal D}(\mathbb{R}^d)$, then (iii)
implies $\widehat{D}^{(n)}(0) = 0$, for all $n = 1,2,\dots$, where $\widehat{D}$ denotes the Fourier transform of
$D$. It follows $\widehat{D}=\widehat{D}(0) = c$ for some constant $c$, since $\widehat{D}$
is an entire function on $\mathbb{C}^d$, by the Paley-Wiener Theorem (Bremermann [12], Theorem 8.28, p. 97). On
the other hand, $D\in {\cal D}(\mathbb{R}^d)\subset S(\mathbb{R}^d)$ implies $\widehat{D}|\mathbb{R}^d\in S(\mathbb{R}^d)$ since
$S(\mathbb{R}^d)$ is closed under Fourier transform. Thus, it follows $c = 0$, i.e. $\widehat{D} = 0$ which implies
$D = 0$ contradicting (ii).

For other classes of nonstandard delta functions we refer to (Robinson [15], p. 133) and to (Todorov [16]).

\smallskip

Our next task is to show the existence of an internal {\it cut-off function}.

\medskip

{\bf Notations}. Let $\Omega$ be an open set of $\mathbb{R}^d$.

1) For any $\varepsilon\in\mathbb{R}_+$, we define
$$
B_\varepsilon= \{x\in \mathbb{R}^d:\| x \| \leq \varepsilon\}\mbox{ and } \Omega_\varepsilon=\{x\in\Omega:d(x,\partial\Omega)\geq\varepsilon\},
$$
where $\|x\|$ is the Euclidean norm in $\mathbb{R}^d$, $\partial\Omega$ is the boundary of $\Omega$ and $d(x,\partial\Omega)$ is the Euclidean distance between $x$ and $\partial\Omega$. We also denote:
$$
{\cal D}_\varepsilon(\Omega) = \{\varphi\in {\cal D}(\Omega): \mbox{supp}\, \varphi\subseteq B_\varepsilon\},
\ \ \ \mathcal E^\prime_\varepsilon(\Omega) = \{T \in \mathcal E^\prime(\Omega) : \mbox{supp}\, T \subseteq \Omega_\varepsilon\}.
$$
2) We shall use the same notation, $\star$, for the convolution operator $\star\colon {\cal D}^\prime(\mathbb{R}^d)\times{\cal D}(\mathbb{R}^d)
\to \mathcal E(\mathbb{R}^d)$ (Vladimirov [13]) and its nonstandard extension
$\star\colon ^*{\cal D}^\prime(\mathbb{R}^d)\times{^*{\cal D}}(\mathbb{R}^d)
\to\, ^* \mathcal E(\mathbb{R}^d)$ as well as for the convolution operator
$\star\colon \mathcal E^\prime_\varepsilon(\Omega)\times {\cal D}_\varepsilon(\Omega)\to   {\cal D}(\Omega)$,
defined for all sufficiently small $\varepsilon\in\mathbb{R}_+$, and for its
nonstandard extension:\, $\star\colon ^* \mathcal E^\prime_\varepsilon(\Omega)\times\, ^* {\cal D}_\varepsilon(\Omega)\to \,^* {\cal D}(\Omega)$,
$\varepsilon\in\, ^*\mathbb{R}_+$, $\varepsilon\approx 0$. We shall use the same notations, $||\cdot||$ and $d(\cdot, \cdot)$, for the Euclidean norm and distance and their nonstandard extensions $^*||\cdot||$ and $^*d(\cdot, \cdot)$, respectively.

3) Let $\tau$ be the usual Euclidean topology on $\mathbb{R}^d$. We denote by $\widetilde{\Omega}$ the set of the nearstandard
points in $^*\Omega$, i.e.
$$
\widetilde{\Omega}= \bigcup_{x\in\Omega}\mu(x), \leqno{(3.2)}
$$
where $\mu(x)$, $x \in \mathbb{R}^d$, is the system of monads of the topological space $(\mathbb{R}^d,\tau)$ (Todorov [6], p. 687).
Recall that if $\xi\in\, {^*\Omega}$, then $\xi\in \widetilde{\Omega}$ if and only if $\xi$ is a finite point whose standard part belongs to
$\Omega$.

\medskip

{\bf Lemma} 3.4. For any positive infinitesimal $\rho$ in $^*\mathbb{R}$ there exists a function $\Pi$ in
$^*{\cal D}(\Omega)$ (a $\rho$-{\it cut-off function}) such that:

\smallskip

a) $\Pi(x)=1$ for all $x\in \widetilde{\Omega}$;

\smallskip

b) $\mbox{supp}\,\Pi\subseteq\Omega_\rho$, where $\Omega_\rho=\{\xi\in\,^*\Omega:\, d(\xi,\partial\Omega)\geq\rho\}$.

\medskip

{\bf Proof.} Let $\rho$ be a positive infinitesimal in $^*\mathbb{R}$ and $D$ be a $\rho$-delta function. Define the internal set
$X= \{\xi\in\,^*\Omega:\, \|\xi\|\leq 1/\rho,\, d(\xi,\partial\Omega)\geq 2\rho\}$ and let $\chi$ be its characteristic function.
Then the function $\Pi=\chi\star D$ has the desired property. $\square$
\section{The Algebra of Asymptotic Functions}

\quad\; We define and study the algebra $^\rho \mathcal E(\Omega)$ of asymptotic functions on an open set $\Omega$ of $\mathbb{R}^d$.
The construction of the algebra $^\rho \mathcal E(\Omega)$, presented here, is a generalization and a refinement of the
constructions in [9] and [10] (by the authors of this paper, respectively), where the algebra $^\rho\mathcal E(\mathbb{R}^d)$ was
introduced by somewhat different but equivalent definitions. On the other hand, the algebra of
asymptotic functions $^\rho\mathcal E(\Omega)$ is somewhat similar to but different from the Colombeau
algebras of new generalized functions (Colombeau [7], [8]). This essential difference between $^\rho\mathcal E(\Omega)$ and Colombeau's
algebras of generalized functions is the properties of the generalized scalars: the scalars of the algebra
$^\rho\mathcal E(\Omega)$ constitutes an algebraically closed field (as any scalars should do) while the scalars of 
Colombeau's algebras are rings with zero divisors (Colombeau [8], \S 2.1). This improvement
compared with J. F. Colombeau's theory is due to the involvement of the nonstandard analysis.

Let $\Omega$ be an open set of $\mathbb{R}^d$ and $\rho\in{^*\mathbb{R}}$ be a positive infinitesimal. We shall keep
$\Omega$ and $\rho$ fixed in what follows. Following (Robinson [2]), we define:

\medskip

{\bf Definition 4.1} (Robinson's Asymptotic Numbers). The field of the complex Robinson $\rho$-asymptotic numbers is defined as the factor space $^\rho\mathbb{C}=\mathbb{C}_\mathcal M/\mathbb{C}_0$, where
\begin{align}
&\mathbb{C}_\mathcal M= \{\xi\in{^*\mathbb{C}}:|\xi|<\rho^{-n}\mbox{ for some } n\in\mathbb{N}\},\notag\\
&\mathbb{C}_0= \{\xi\in{^*\mathbb{C}}:|\xi|<\rho^{n}\mbox{ for all } n\in\mathbb{N}\},
\end{align}
($\mathcal M$ stands for {\it moderate}). We define the embedding $\mathbb{C}\subset{^\rho\mathbb{C}}$ by $c\to q(c)$, where
$q\colon \mathbb{C}_\mathcal M\to {^\rho\mathbb{C}}$ is the quotient mapping.
The field of the real asymptotic numbers is defined by $^\rho\mathbb{R}=q(^*\mathbb{R}\cap\mathbb{C}_\mathcal M)$.

It is easy to check that $\mathbb{C}_0$ is a maximal ideal in $\mathbb{C}_\mathcal M$ and hence $^\rho\mathbb{C}$ is a field.
Also $^\rho\mathbb{R}$ is a real closed totally ordered nonarchimedean field (since $^*\mathbb{R}$ is a real closed totally ordered field) containing $\mathbb{R}$ as a totally ordered subfield. Thus, it follows that $^\rho\mathbb{C}={^\rho\mathbb{R}}(i)$ is an algebraically closed field, where $i=\sqrt{-1}$.

	The algebra of asymptotic functions is, in a sense, a $\mathcal C^\infty$-counterpart of A. Robinson's asymptotic numbers $^\rho\mathbb C$:
	
\medskip

{\bf Definition 4.2} (Asymptotic Functions on $\Omega$). (i) We define the class $^\rho\mathcal E(\Omega)$ of the $\rho$-{\it asymptotic
functions on} $\Omega$ (or simply, {\it asymptotic functions on} $\Omega$ if no confusion could arise) as the factor space
$^\rho\mathcal E(\Omega)=\mathcal E_\mathcal M(\Omega)/\mathcal E_0(\Omega)$, where
$$
\mathcal E_\mathcal M(\Omega)=\{f\in\, {^*\mathcal E(\Omega)}: \partial^\alpha f(\xi)\in \mathbb{C}_\mathcal M \mbox{ for all } \alpha\in \mathbb{N}_0^d\mbox{ and all } \xi\in\widetilde{\Omega}\},
$$
$$
\mathcal E_0(\Omega)=\{f\in{^*\mathcal E(\Omega)}: \partial^\alpha f(\xi)\in \mathbb{C}_0 \mbox{ for all } \alpha\in \mathbb{N}_0^d\mbox{ and all } \xi\in\widetilde{\Omega}\},
$$
and $\widetilde{\Omega}$ is the set of the nearstandard points of $^*\Omega$ (3.2). The functions in $\mathcal E_\mathcal M(\Omega)$ are called
$\rho$-{\it moderate} (or, simply, {\it moderate}) and those in $\mathcal E_0(\Omega)$ are called $\rho$-{\it null functions} (or, simply,
{\it null functions}).

(ii) The pairing between $^\rho{\mathcal E}(\Omega)$ and ${\mathcal D}(\Omega)$ with values in $^\rho\mathbb{C}$, is defined by
$$
\langle Q_\Omega(f),\varphi\rangle=q\left(\int_{^*\Omega}f(x)^*\!\varphi(x)\, dx\right),
$$
where $q\colon \mathbb{C}_\mathcal  M\to\, ^\rho\mathbb{C}$ and $Q_\Omega\colon \mathcal E_\mathcal  M(\Omega)\to\,^\rho\mathcal E(\Omega)$ are the corresponding quotient mappings,
$\varphi$ is in ${\cal D}(\Omega)$ and $^*\!\varphi$ is its nonstandard extension of $\varphi$.

(iii) We define the {\it canonical embedding} $\mathcal E(\Omega)\subset{^\rho \mathcal E(\Omega)}$ the mapping
$\sigma_\Omega\colon f\to Q_\Omega(^*\! f)$, where $^*\!f$
is the nonstandard extension of $f$.

\medskip

{\bf Example 4.3.} Let $D$ be a nonstandard $\rho$-delta function in the sense of Definition 3.2. Then
$D\in\mathcal E_\mathcal M(\mathbb{R}^d)$. In addition, $D|{^*\Omega\in\mathcal E_\mathcal M}(\Omega)$, where $D|^*\Omega$ denotes the pointwise restriction of $D$ on $^*\Omega$. To show this, denote $\displaystyle|\ln\rho|^{-1}\left(\rho^{d+|\alpha|}\sup_{x\in{^*\mathbb{R}^d}}|\partial^\alpha D(x)| \right)=h_\alpha$
and observe that $h_\alpha\approx 0$ for all $\alpha\in\mathbb{N}_0^d$ by the definition of $D$.
Thus, for any (finite) $x$ in $^*\mathbb{R}^d$ and any $\alpha\in\mathbb{N}_0^d$ we have
$|\partial^\alpha D(x)|\leq\sup_{x\in{^*\mathbb{R}}}|\partial^\alpha D(x)|=\frac{h_\alpha|\ln\rho|}{\rho^{d+|\alpha|}}<\rho^{-n}$, for $n=d+|\alpha|+1$, thus, $D\in\mathcal E_\mathcal M(\mathbb{R}^d)$. On the other hand, $D|{^*\Omega}\in\mathcal E_\mathcal  M(\Omega)$ follows immediately from the fact that
$\widetilde{\Omega}$ consists of finite points in $^*\mathbb{R}^d$ only.

\medskip

{\bf Theorem 4.4} (Differential Algebra). (i) The class of asymptotic functions $^\rho \mathcal E(\Omega)$ is a differential
algebra over the field of the complex asymptotic numbers $^\rho\mathbb{C}$.

(ii) $\mathcal E(\Omega)$ is a {\it differential subalgebra} of $^\rho \mathcal E(\Omega)$ over the scalars $\mathbb{C}$ under the canonical embedding
$\sigma_\Omega$. In addition, $\sigma_\Omega$ preserves the pairing in the sense that $\langle f,\varphi\rangle = \langle \sigma_\Omega(f),\varphi\rangle$ for all $f$ in $\mathcal E(\Omega)$ and for
all $\varphi$ in ${\cal D}(\Omega)$, where
$\langle f,\varphi\rangle = \int_{\Omega}f(x)\varphi(x)\,dx$  is the usual pairing between $\mathcal E(\Omega)$ and ${\cal D}(\Omega)$.

\medskip

{\bf Proof.} (i) It is clear that $\mathcal E_\mathcal M(\Omega)$ is a differential ring and $\mathcal E_0(\Omega)$ is a differential ideal in $\mathcal E_\mathcal M(\Omega)$
since $\mathbb{C}_\mathcal M$ is a ring and  $\mathbb{C}_0$ is an ideal in  $\mathbb{C}_\mathcal M$ and, on the other hand, both $\mathcal E_\mathcal M(\Omega)$ and
$\mathcal E_0(\Omega)$ are closed under differentiation, by definition. Hence, the factor space $^\rho\mathcal E(\Omega)$ is also a differential ring. It is clear that, $\mathcal E_\mathcal M(\Omega)$ is a module over the ring $\mathbb{C}_\mathcal M$, and, in addition, the annihilator $\{c\in \mathbb{C}_\mathcal M: (\forall f\in \mathcal E_\mathcal M(\Omega))(cf\in\mathcal E_0(\Omega))\}$
of $\mathbb{C}_\mathcal M$ coincides with the ideal $\mathbb{C}_0$. Thus, $^\rho \mathcal E(\Omega)$ becomes an algebra over the field of the complex
asymptotic numbers $^\rho\mathbb{C}$.

(ii) Assume that $\sigma_\Omega(^*\! f)=0$ in $^\rho \mathcal E(\Omega)$, i.e. $^*\! f\in\mathcal E_0(\Omega)$.
By the definition of $\mathcal E_0(\Omega)$ (applied for $\alpha=0$ and $n = 1$), it follows $f =0$ since $^*\! f$ is an extension of $f$ and
$\rho$ is an infinitesimal. Thus, the mapping $f\to\sigma_\Omega(f)$ is injective. It preserves the algebraic operations since the mapping
$f\to\, ^*\! f$ preserves them. The preserving of the pairing follows immediately from the fact that $\int_{^\star\Omega}\, ^*\! f(x)\, dx=\int_{\Omega}f(x)\, dx$, by the Transfer Principle (Todorov [6], p. 686). The proof is complete. $\square$


\section{Embedding of Schwartz Distributions}

\quad\; Let $\Omega$ be (as before) an open set of $\mathbb{R}^d$. Recall that the {\it Schwartz embedding}
$L_\Omega\colon {\cal L}_{loc}(\Omega)\to {\cal D}^\prime(\Omega)$
from ${\cal L}_{loc}(\Omega)$ into ${\cal D}^\prime(\Omega)$ is defined by the formula:
$$
\langle L_\Omega(f),\varphi\rangle=\int_{\Omega}f(x)\varphi(x)\, dx, \ \ \ \varphi\in{\cal D}(\Omega). \leqno{(5.1)}
$$
Here ${\cal L}_{loc}(\Omega)$ denotes, as usual, the space of the locally (Lebesgue) integrable complex valued functions
on $\Omega$ (Vladimirov [13]). The Schwartz embedding $L_\Omega$ preserves the addition and multiplication by a
complex number, hence, the space ${\cal L}_{loc}(\Omega)$ can be considered as a linear subspace of ${\cal D}^\prime(\Omega)$.
In addition, the restriction $L_\Omega | \mathcal E(\Omega)$ of $L_\Omega$ on $\mathcal E(\Omega)$ (often denoted also by $L_\Omega$) preserves the partial differentiation of any order and in this sense $\mathcal E(\Omega)$ is a differential linear subspace of ${\cal D}^\prime(\Omega)$.
In short, we have the chain of linear embeddings: $\mathcal E(\Omega)\subset{\cal L}_{loc}(\Omega)\subset {\cal D}^\prime(\Omega)$.

The purpose of this section is to show that the algebra of asymptotic functions $^\rho \mathcal E(\Omega)$ contains an
isomorphic copy of the space of Schwartz distributions ${\cal D}^\prime(\Omega)$ and, hence, to offer {\it a solution of the
Problem of Multiplication of Schwartz Distributions}. This result is a generalization of some results in [9]
and [10] (by the authors of this paper, respectively) where only the embedding of the tempered
distributions $S^\prime(\mathbb{R}^d)$ in $^\rho\mathcal E(\mathbb{R}^d)$ has been established. The embedding of all distributions
${\cal D}^\prime(\Omega)$ discussed here, presents an essentially different situation.

	The spaces $\widetilde{\mathcal E}(\Omega)$ and $\widetilde{{\cal D}}(\Omega)$, defined below, are immediate generalizations of the spaces
$\widetilde{\mathcal E}(\mathbb{R}^d)$ and $\widetilde{{\cal D}}(\mathbb{R}^d)$, introduced in (Stroyan and Luxemburg [17], (10.4), p. 299):
$$
\widetilde{\mathcal E}(\Omega) = \{\varphi\in{^*\mathcal E}(\Omega): \partial^\alpha\varphi(x) \mbox{ is a finite number in } {^*\mathbb{C}}
\mbox{ for all } x\in\widetilde{\Omega}\mbox{ and all } \alpha\in\mathbb{N}_0^d\}, \leqno{(5.2)}
$$
$$
\widetilde{{\cal D}}(\Omega) = \{\varphi\in{^*\mathcal E(\Omega)}: \partial^\alpha\varphi(x) \mbox{ is a finite number in } ^*\mathbb{C}
\mbox{ for all }
$$
$$x\in\widetilde{\Omega}, \alpha\in\mathbb{N}_0^d \mbox{ and } \varphi(x)=0 \mbox{ for all } x\in{^*\Omega}\backslash\widetilde{\Omega} \},\leqno{(5.3)}
$$
Obviously, we have $\widetilde{{\cal D}}(\Omega)\subset\widetilde{{\mathcal E}}(\Omega)\subset\mathcal E_\mathcal M(\Omega)$.
Notice as well that $\varphi\in \widetilde{{\cal D}}(\Omega)$ implies $\varphi\in{^*{\cal D}}(G)$ for
some open relatively compact set $G$ of $\Omega$. We have also the following simple result:

\medskip

{\bf Lemma 5.1.} If $T\in {\cal D}^\prime(\Omega)$ and $\varphi\in \mathcal E_0(\Omega)\cap\, \widetilde{{\cal D}}(G)$, then
$\langle {^* T},\varphi\rangle\in\mathbb{C}_0$.

{\bf Proof.} Observe that $\mathcal E_0(\Omega)\cap \widetilde{{\cal D}}(\Omega)$ implies $\varphi\in \mathcal E_0(\Omega)\cap{^*{\cal D}(G)}$ for some open relatively
compact set $G$ of $\Omega$. By the continuity of $T$ (and Transfer Principle) there exist constants $C\in\mathbb{R}_+$, and
$m \in\mathbb{N}_0$ such that
$$
|\langle {^*T},\varphi\rangle |\leq C\sum_{|\mu|\leq m}\sup_{x\in{^* G}}|\partial^\mu\varphi(x)|
$$
On the other hand, $\displaystyle C\sum_{|\mu\leq m}\sup_{x\in^\star G}|\partial^\mu\varphi(x)| < \rho^n$ for all $n\in\mathbb{N}$, since
$\varphi\in\mathcal E_0(\Omega)$, by assumption. Thus,
$|\langle {^* T},\varphi\rangle |<\rho^n$ for all $n\in\mathbb{N}$. $\square$

Let $D$ be a $\rho$-delta function in the sense of Definition 3.2. We shall keep $D$ (along with $\Omega$ and $\rho$)
fixed in what follows.

\medskip

{\bf Definition 5.2} ({\it Embedding of Schwartz Distributions}). We define the embedding
${\cal D}^\prime(\Omega)\subset\, ^\rho \mathcal E(\Omega)$ by $\Sigma_{D,\Omega}\colon T\to Q_\Omega((^* T\,\Pi_\Omega)\star D)$, where
$^* T$ is the nonstandard extension of $T$, $\Pi_\Omega$ is
a (an arbitrarily chosen) $\rho$-cut-off function for $\Omega$ (Lemma 3.4), $^* T\,\Pi_\Omega$ is the Schwartz product between
$^*T$ and $\Pi_\Omega$ in $^*{\cal D}^\prime(\Omega)$ (applied by Transfer Principle), $\star$ is the convolution operator and
$Q_\Omega\colon\mathcal E_\mathcal M(\Omega)\to \,^\rho\mathcal E(\Omega)$ is the quotient mapping in the definition of $^\rho\mathcal E(\Omega)$ (Definition 4.2).

The cut-off function $\Pi_\Omega$ can be dropped in the above definition, i.e. $\Sigma_{D,\Omega}\colon T\to Q_\Omega({^* T}\star D)$, in
some particular cases; e.g. when:

a) $T$ has a compact support in $\Omega$;

b) $\Omega=\mathbb{R}^d$.

{\bf Proposition 5.3} (Correctness). $T \in {\cal D}^\prime(\Omega)$ implies $({^* T}\,\Pi_\Omega)\star D\in \mathcal E_\mathcal M(\Omega)$.

{\bf Proof.} Choose $\alpha\in\mathbb{N}_0^d$ and $x\in\widetilde{\Omega}$.
Since we have $\partial^\alpha (({^*T}\, \Pi_\Omega)\star D)(x)=(\partial^\alpha(^* T)\star D)(x)$
(by the definition of $\Pi_\Omega$), we need to show that $\partial^\alpha({^*T}\star D)(x)\in\mathbb{C}_\mathcal M$ only, i.e. that
$|\partial^\alpha(^*T)\star D)(x)|<p^{-m}$ for some $m\in\mathbb{N}$ ($m$ might depend on $\alpha$). We start with the case
$\alpha = 0$. Denote $D_x(\xi)=D(\xi-x)$, $\xi\in{^*\mathbb{R}}$ and observe that $\mbox{supp}(D_x)\subseteq{^*G}$ for some open relatively compact set $G$ of $\Omega$, since $D_x$, vanishes on $^*\Omega\setminus\widetilde{\Omega}$.
Next, by the continuity of $T$ (and the Transfer Principle), there
exist constants $m\in\mathbb{N}_0$ and $C\in\mathbb{R}_+$, such that
$$
|({^*T}\star D)(x)| = |\langle ^* T, D_x|^*\Omega\rangle|\leq C\sum_{|\mu|\leq m}\sup_{\xi\in{^* G}}|\partial_\xi^\mu D(x-\xi)|
$$
Finally, there exists $n \in N$ such that $\displaystyle\sum_{|\mu|\leq m}\sup_{\xi\in^\star G}|\partial_\xi^\mu D(x-\xi)|<p^{-n}$, since $D|^\star G$ is a $\rho$-moderate function (Example 4.3). Combining these arguments, we have: $|(^*T\star D)(x)|\leq C\rho^{-n}<\rho^{-(n+1)}$, as
required. The generalization for arbitrary multiindex $\alpha$ follows immediately since
$\partial^\alpha({^*T}\star D)=(\partial^\alpha(^*T))\star D)={^* (\partial^\alpha} T)\star D$, by Transfer Principle, and
$\partial^\alpha T$ is (also) in ${\cal D}^\prime(\Omega)$. $\square$

\medskip

{\bf Proposition 5.4.} $f \in \widetilde{\mathcal E}(\Omega)$ implies $(f\Pi_\Omega)\star D-f\in \mathcal E_0(\Omega)$.

\smallskip

{\bf Proof.} Let $x\in\widetilde{\Omega}$ and $\alpha\in\mathbb{N}_0^d$. Since we have $\partial^\alpha[((f\Pi_\Omega)\star D)(x)-f(x)]=
\partial^\alpha[(f\star D)(x)-f(x)]$ (by the definition of $\Pi_\Omega$), we need to show that $\partial^\alpha[(f\star D)(x)-f(x)]\in\mathbb C_0$
only. Choose $n \in\mathbb{N}$. We need to show that $|\partial^\alpha[(f\star D)(x)-f(x)]|<\rho^n$. We start first with the case
$\alpha=0$. By Taylor's formula (applied by transfer), we have
$$
f(x-\xi)-f(x)=\sum_{|\beta|=1}^{n}\frac{(-1)^{|\beta|}\partial^\beta f(x)}{\beta!}\,\xi^\beta+\frac{(-1)^{n+1}}{(n+1)!}\sum_{|\beta|=n+1}\partial^\beta f(\eta(\xi))\,\xi^\beta
$$
for any $\xi\in\widetilde{\Omega}$, where $\eta(\xi)$ is a point in $^*\Omega$ “between $x$ and $\xi$.".
Notice that the point  $\eta(\xi)$ is also in $\widetilde{\Omega}$. It follows
$$
(f\star D)(x)-f(x)=\int_{\|\xi\|\leq \rho}D(\xi)[f(x-\xi)-f(x)]\, d\xi=\frac{(-1)^{n+1}}{(n+1)!}\sum_{|\beta|=n+1}\int_{\|\xi\|\leq \rho}
D(\xi)\xi^\beta\partial^\beta f(\eta(\xi))\,d\xi
$$
since $\int_{\|\xi\|\leq \rho}D(\xi)\xi^\beta \, d\xi=0$, by the definition of $D$. Thus, we have
$$
|(f\star D)(x)-f(x)|\leq \frac{\rho^{n+1}}{(n+1)!}\left(\int_{^\star\mathbb{R}^d}|D(x)|\,dx\right)\left(\sum_{|\beta|=n+1}\sup_{\|\xi\|\leq\rho}
|\partial^\beta f(\eta(\xi)) |\right)<\rho^n,
$$
as desired, since, on one hand, $\int_{^\star\mathbb{R}^d}|D(x)|\, dx\approx 1$, by the definition of $D$ and on the other hand, the
above sum is a finite number because $^*\partial^\beta f(\eta(\xi))$ are all finite due to $\eta(\xi)\in \widetilde{\Omega}$.
The generalization for an arbitrary $\alpha$ is immediate since
$\partial^\alpha[(f\star D)(x)-f(x)]=(\partial^\alpha f\star D)(x)-\partial^\alpha f(x)$, by the Transfer
Principle. $\square$

\medskip

{\bf Corollary 5.5.} (i) $f \in\mathcal E(\Omega)$ implies $(^* f\,\Pi_\Omega)\star D-\, ^*\! f\in\mathcal E_0(\Omega)$.

(ii) $\varphi \in  {\cal D}(\Omega)$ implies $({^*\!\varphi}\, \Pi_\Omega)\star D-\,{^*\!\varphi}\in \mathcal E_0(\Omega)\cap\widetilde{{\cal D}}(\Omega)$.

{\bf Proof.} (i) follows immediately from the above proposition since $f \in\mathcal E(\Omega)$ implies $^*\! f \in \widetilde{\mathcal E}(\Omega)$.
(ii) Both $(^*\!\varphi\,\Pi_\Omega)\star D$ and $^*\!\varphi$ vanish on $^*\Omega\setminus\widetilde{\Omega}$ since their supports are within an open relatively compact neighborhood $G$ of $\mbox{supp}(\varphi)$ and the latter is a compact set of $\Omega$, by assumption. Thus, $(^*\!\varphi\,\Pi_\Omega)\star D-{^*\!\varphi}\in{^*\mathcal D(G)}\subset\tilde{\mathcal D}(\Omega)$ as required. $\square$

	\quad\; Denote $\check{D}(x) = D(- x)$ and recall that $\check{D} = D$ since $D$ is symmetric (Definition 3.2).

\medskip

{\bf Proposition 5.6.} If $T \in {\cal D}^\prime(\Omega)$ and $\varphi\in {\cal D}(\Omega)$, then
$$
\int_{^*\Omega}(({^* T}\,\Pi_\Omega)\star D)(x)\, {^*\!\varphi}(x)\, dx-\langle T,\varphi\rangle\in \mathbb{C}_0
$$
{\bf Proof.} Using the properties of the convolution operator (applied by transfer), we have
$$
\int_{^*\Omega}(({^* T}\,\Pi_\Omega)\star D)(x)\,{^*\!}\varphi(x)\, dx-\langle T,\varphi\rangle
$$
$$
=\langle (^*T\, \Pi_\Omega)\star D,{^*\!\varphi}\, \rangle-\langle {^* T},{^*\!\varphi}\rangle=
\langle {^* T}\, \Pi_\Omega, {^*\! \varphi}\star \check{D}\rangle-\langle {^* T}\,\Pi_\Omega,{^*\!\varphi}\rangle
$$
$$
=\langle {^* T}\, \Pi_\Omega,\,{^*\!\varphi}\star \check{D}-\,{^*\!\varphi}\rangle=\langle ^*T,\,{^*\!\varphi}\star D-\,{^*\!\varphi}\rangle\in \mathbb{C}_0,
$$
by Lemma 5.1 since $^*\!\varphi\star D -\,{^*\!\varphi} \in  \mathcal E_0(\Omega)\cap \widetilde{\mathcal D}(\Omega)$, by Corollary 5.5. $\square$

We are ready to state our main result:

\medskip

{\bf Theorem 5.7} (Properties of $\Sigma_{D,\Omega}$). (i) $\Sigma_{D,\Omega}$ preserves the pairing in the sense that for all $T$ in
${\cal D}^\prime(\Omega)$ and all $\varphi$ in ${\cal D}(\Omega)$ we have $\langle T,\varphi\rangle =\langle \Sigma_{D,\Omega}(T),\varphi\rangle$, where the left hand side is the (usual) pairing
of $T$ and $\varphi$ in ${\cal D}^\prime(\Omega)$, while the right hand side is the pairing of $\Sigma_{D,\Omega}(T)$ and
$\varphi$ in $^\rho\mathcal E(\Omega)$ (Definition 4.2).

(ii) $\Sigma_{D,\Omega}$ is {\it injective} and it preservers {\it all linear operations} in ${\cal D}^\prime(\Omega)$: the addition, multiplication by (standard) complex numbers and the partial differentiation of any (standard) order.

(iii) $\Sigma_{D,\Omega}$ is {\it an extension of the canonical embedding} $\sigma_\Omega$ defined earlier in Definition 4.2 in the
sense that $\sigma_{\Omega}=\Sigma_{D,\Omega}\circ L_\Omega$, where $L_\Omega$ is the Schwartz embedding (5.1) restricted on
$\mathcal E(\Omega)$ and $\circ$ denotes composition. Or, equivalently, the following diagram is commutative: 
$$
\begin{tikzcd}[column sep=small]
\mathcal E(\Omega) \arrow{r}{L_\Omega}  \arrow{rd}[swap]{\sigma_\Omega} 
  & \mathcal D^\prime(\Omega) \arrow{d}{\Sigma_{D, \Omega}} \\
    & {^\rho\mathcal E(\Omega)}.
\end{tikzcd}
$$
{\bf Proof.} (i) Denote (as before) $\check{D}(x) = D(-x)$ and recall that $\check{D} = D$ (Definition 3.2). We

have
\begin{eqnarray*}
\langle\Sigma_{D,\Omega}(T),\varphi\rangle&=&\langle Q_\Omega((^* T\, \Pi_\Omega)\star D),\varphi\rangle-\langle T,\varphi\rangle\\
&=&q\left(\int_{^*\Omega}((\Pi_\Omega\, {^*T})\star D)(x)\,{^*\!\varphi}(x)\,dx \right)-q (\langle T,\varphi\rangle)\\
&=&q\left(\int_{^*\Omega}((\Pi_\Omega\, {^* T})\star D)(x)\, {^*\!\varphi}(x)-\langle T,\varphi\rangle\right)=0,
\end{eqnarray*}
because $\int_{^*\Omega}((^* T\, \Pi_\Omega)\star D)(x)\,{^*\!\varphi}(x)\, dx - \langle T,\varphi\rangle \in \mathbb{C}_0$, by Proposition 5.6. Here $\langle T,\varphi\rangle=q(\langle T,\varphi\rangle)$
holds because $\langle T,\varphi\rangle$ is a standard (complex) number.

(ii) The injectivity of $\Sigma_{D,\Omega}$ follows from (i):
$\Sigma_{D,\Omega}(T)=0$ in $^\rho \mathcal E(\Omega)$ implies $\langle\Sigma_{D,\Omega}(T),\varphi\rangle = 0$ for
all $\varphi \in {\cal D}(\Omega)$, which is equivalent to $\langle T,\varphi\rangle = 0$ for all $\varphi\in {\cal D}(\Omega)$, by (i), thus,
$T=0$ in ${\cal D}^\prime(\Omega)$, as required. The preserving of the linear operations follows from the fact that both the extension mapping
$^*$ and the convolution $\star$  (applied by Transfer Principle) are linear operators.

(iii) For any $f\in \mathcal E(\Omega)$ we have $\sigma(f)=Q_{\Omega}(^*\! f)=Q_{\Omega}((^*\! f\, \Pi_{\Omega})\star D)=
Q_{\Omega}((^*L_\Omega(f)\, \Pi_{\Omega})\star D)=\Sigma_{D,\Omega}(L_\Omega(f))$, as required,
since $^*\!f-({^*\! f}\, \Pi_{\Omega})\star D\in \mathcal E_0(\Omega)$, by Corollary 5.5. $\square$

\medskip

{\bf Remark 5.8} (Multiplication of Distributions). As a consequence of the above result, the
Schwartz distributions in ${\cal D}^\prime(\Omega)$ can be multiplied within the associative and commutative differential
algebra $^\rho \mathcal E(\Omega)$ (something impossible in ${\cal D}^\prime(\Omega)$ itself). By the property (iii) above, the multiplication in
$^\rho \mathcal E(\Omega)$ coincides on $\mathcal E(\Omega)$ with the usual (pointwise) multiplication in $\mathcal E(\Omega)$. Thus, the class $^\rho \mathcal E(\Omega)$, endowed with an embedding $\Sigma_{D,\Omega}$, presents {\it a solution of the problem of multiplication of Schwartz
distributions} which, in a sense, is optimal, in view of the Schwartz impossibility results (Schwartz [1])
(for a discussion we refer also to (Colombeau [7], \S 2.4) and (Oberguggenberg [18], \S2). We
should mention that the existence of an embedding of ${\cal D}^\prime(\mathbb{R}^d)$ into $^\rho\mathcal E(\mathbb{R}^d)$ can be proved also by sheaf-theoretical arguments as indicated in (Oberguggenberger [18], \S 23).

\medskip

{\bf Remark 5.9} (Nonstandard Asymptotic Analysis). We sometimes refer to the area connected
directly or indirectly with the fields $^\rho\mathbb{R}$ as Nonstandard Asymptotic Analysis. The fields $^\rho\mathbb{R}$ were
introduced in Robinson [2] and are sometimes known as ``Robinson's nonarchimedean valuation
fields''. The terminology ``Robinson's asymptotic numbers'',  chosen in this paper, is due to the role of $^\rho\mathbb{R}$
for the asymptotic expansions of classical functions (Lightstone and Robinson [3]) and also to
stress the fact that in our approach $^\rho\mathbb{C}$ plays the role of the scalars of the algebra $^\rho\mathcal E(\Omega)$. Linear spaces
over the field $^\rho\mathbb{R}$ has been studied in Luxemburg [19] in order to establish a connection between
nonstandard and nonarchimedean analysis. More recently $^\rho\mathbb{R}$ has been used in Pestov [20] for
studying Banach spaces. The field $^\rho\mathbb{R}$ has been exploited by Li Bang-He [21] for multiplication of
Schwartz distributions.


\end{document}